\documentclass[12pt,reqno]{amsart}

\usepackage{amssymb}

\usepackage{mathrsfs}

\DeclareMathAlphabet{\mathpzc}{OT1}{pzc}{m}{it}

\usepackage[all]{xy}
\entrymodifiers={+!!<0pt,\fontdimen22\textfont2>}
\font\sss=cmss8

\def\fa{{\mathfrak a}}

\def\frakd{{\mathfrak d}}

\def\sC{\mbox{\sf C}}
\def\sD{\mbox{\sf D}}

\def\sK{\mbox{\sf K}}

\def\sT{\mbox{\sf T}}

\def\D{\mbox{\sD}}

\def\Dsmall{\mbox{\sss D}}
\def\dual{\operatorname{D}}

\def\Ext{\operatorname{Ext}}

\def\finpres{\mbox{\sf fp}}

\def\gldim{\operatorname{gldim}}

\def\H{\operatorname{H}}

\def\Hom{\operatorname{Hom}}
\def\id{\operatorname{id}}
\def\Image{\operatorname{Im}}

\def\Inj{\mbox{\sf Inj}}

\def\Ker{\operatorname{Ker}}

\def\mod{\mbox{\sf mod}}

\def\prod{\operatorname{prod}}

\numberwithin{equation}{part}

\newtheorem{Lemma}{Lemma}[section]
\newtheorem{Theorem}[Lemma]{Theorem}

\theoremstyle{definition}
\newtheorem{Definition}[Lemma]{Definition}
\newtheorem{Setup}[Lemma]{Setup}

\newtheorem{Remark}[Lemma]{Remark}

\def\cat{\sT}
\def\catsmall{\mbox{\sss T}}
\def\subcat{\sC}

\begin{document}

\setlength{\parindent}{0pt}
\setlength{\parskip}{7pt}
%The default \baselineskip is close to 4.8mm
%\setlength{\baselineskip}{5.3mm}

\title[Auslander-Reiten triangles in subcategories]
{Auslander-Reiten triangles in subcategories}

\author{\ \ Peter J\o rgensen}
\address{Department of Pure Mathematics, University of Leeds,
Leeds LS2 9JT, United Kingdom}
\email{popjoerg@maths.leeds.ac.uk}
\urladdr{http://www.maths.leeds.ac.uk/\~{ }popjoerg}
%\thanks{Date: \today. A thank you would go here}

%\author{Next author goes here}
%\address{Next author's postal address goes here}
%\email{Next author's mail address goes here}


%14A22 = Noncommutative algebraic geometry
%16D90 = Module categories [See also 16Gxx, 16S90]; module theory in a category-theoretic context; Morita equivalence and duality
%16E05 = Syzygies, resolutions, complexes
%16E30 = Homological functors on modules (Tor, Ext, etc.)
%16E45 = Differential graded algebras and applications
%16E65 = Homological conditions on rings (generalizations of regular,
%        Gorenstein, Cohen-Macaulay rings, etc.)
%16G70 = Auslander-Reiten sequences
%16W50 = Graded rings and modules
%18E30 = Derived categories, triangulated categories
%55P62 = Rational homotopy theory

\begin{abstract} 

This paper introduces Auslander-Reiten triangles in subcategories of
triangulated categories.

The main theorem of the paper shows that there is a close connection
with covers and envelopes, also known as minimal right- and
left-ap\-proxi\-ma\-ti\-ons.  Namely, under suitable assumptions, if
$C$ is an object in the sub\-ca\-te\-go\-ry $\subcat$ of the
triangulated category $\cat$ and
\[
  \xymatrix{
    X \ar[r] & Y \ar[r] & C \ar[r] & {}
           }
\]
is an Auslander-Reiten triangle in $\cat$, then there is an
Auslander-Reiten triangle
\[
  \xymatrix{
    A \ar[r] & B \ar[r] & C \ar[r] & {}
           }
\]
in $\subcat$ if and only if there is a $\subcat$-cover of the form
$\xymatrix{ A \ar[r] & X. }$

The main theorem is used to give a new proof of the existence of
Auslander-Reiten sequences over finite dimensional algebras.

\end{abstract}

\maketitle

% \setcounter{section}{-1}
% \section{Introduction}
% \label{sec:introduction}

\section{Introduction}
\label{sec:introduction}

\noindent
Auslander-Reiten sequences, also known as almost split exact
se\-quen\-ces, are one of the central tools in the representation
theory of Artin algebras.  The precise definition can be found in many
places; the book \cite{ARS} gives an excellent presentation.  For the
purposes of this introduction, the best way to treat these sequences
is perhaps to give a sample of their properties.

If
\begin{equation}
\label{equ:m}
  \xymatrix{
    0 \ar[r] & A \ar[r] & B \ar[r] & C \ar[r] & 0
           }
\end{equation} 
is an Auslander-Reiten sequence in the category of finitely generated
modules over an Artin algebra, then $A$ and $C$ are indecomposable
mo\-du\-les which determine each other.  Conversely, if $A$ or $C$ is
a given indecomposable module, then there is an Auslander-Reiten
sequence \eqref{equ:m}.  Moreover, the homomorphisms in the sequence
determine the so-called irreducible homomorphisms out of $A$ and into
$C$, and in good cases, the irreducible homomorphisms between
indecomposable mo\-du\-les can be used to determine all module
homomorphisms.  In effect, the Auslander-Reiten sequences sum up a
great deal of information about the module category of an Artin
algebra.

The theory of Auslander-Reiten sequences was founded in
\cite{AusArtIII}, and in fact, makes sense for abstract abelian
categories.  Later on, Auslander and Smal\o\ in \cite{AuslanderSmalo}
developed a more general theory of Auslander-Reiten sequences in
subcategories; these permit the study of non-abelian subcategories of
abelian categories.

Another later advance was Happel's theory of Aus\-lan\-der-Rei\-ten
triangles, see \cite{Happel}; these play the same role in triangulated
categories as Auslander-Reiten sequences do in abelian categories.

On this background, it seems very natural to try to develop a theory
of Auslander-Reiten triangles in subcategories.  Such a theory is
given in this paper; it permits the study of non-tri\-an\-gu\-la\-ted
subcategories of triangulated ca\-te\-go\-ri\-es.

%\begin{center}
%\begin{tabular}{l|cc}
%\rule[-3ex]{0ex}{7ex}                   &  \begin{minipage}{23ex}
%        \centerline{Abelian}
%        \centerline{categories}
%      \end{minipage}                      &  \begin{minipage}{12ex}
%        \centerline{Triangulated}
%        \centerline{categories}
%      \end{minipage} \\ \cline{1-3}
%\rule{0ex}{5ex} \begin{minipage}{25ex}
%        \leftline{Auslander-Reiten}
%        \leftline{theory}
%      \end{minipage}  & Auslander-Reiten \cite{AusArtIII}       & Happel \cite{Happel} \\
%\rule{0ex}{5ex} \begin{minipage}{25ex}
%        \leftline{Auslander-Reiten}
%        \leftline{theory for subcategories}
%      \end{minipage}    & Auslander-Smal\o\ \cite{AuslanderSmalo} & \\
%\end{tabular}
%\end{center}

Auslander-Reiten triangles in subcategories are introduced in
Definition \ref{def:arsubcat}.  The main theorem I will prove about
them is Theorem \ref{thm:main}, which shows that there is a close
connection between Auslander-Reiten triangles in a subcategory and the
approximation properties of the subcategory.  Namely, let $\cat$ be a
suitable triangulated category, $\subcat$ a suitable subcategory.  If
$C$ is in $\subcat$ and there is an Auslander-Reiten triangle
\[
  \xymatrix{
    X \ar[r] & Y \ar[r] & C \ar[r] & {}
           }
\]
in $\cat$, then there is an Auslander-Reiten triangle
\[
  \xymatrix{
    A \ar[r] & B \ar[r] & C \ar[r] & {}
           }
\]
in $\subcat$ if and only if there is a $\subcat$-cover, that is, a
minimal right-$\subcat$-approximation, of the form $\xymatrix{ A
\ar[r] & X; }$ see Definition \ref{def:cover_and_envelope}.  This
echoes the properties of Auslander-Reiten sequences in subcategories
dis\-co\-ve\-red by Kleiner in \cite{Kleiner}.

The paper is organized as follows: This introduction ends with some
blanket items, including the definition of Auslander-Reiten triangles
in subcategories.  Section \ref{sec:lemmas} gives some lemmas.
Section \ref{sec:main} proves the main theorem, stating the connection
between Auslander-Reiten triangles in subcategories and the
approximation properties of the subcategory.

Section \ref{sec:example} uses the main theorem to give a new proof of
the existence of Auslander-Reiten sequences over finite dimensional
algebras.

The following setup takes place in the world of triangulated
categories.  For background information, see \cite{Happel}
or \cite{Neemanbook}. 

\begin{Setup}
\label{set:blanket}
Throughout, $k$ is a field and $\cat$ is a skeletally small $k$-linear
triangulated category with split idempotents in which each $\Hom$
space is finite dimensional over $k$.

By $\subcat$ is denoted a full subcategory of $\cat$ which is closed
under extensions and direct summands.  That is, if $\xymatrix{ A
\ar[r] & B \ar[r] & C \ar[r] & {} }$ is a distinguished triangle with
$A$ and $C$ in $\subcat$ then $B$ is isomorphic to an object in
$\subcat$, and if $A$ is an object in $\subcat$ for which $A \cong A_1
\amalg A_2$, then $A_1$ and $A_2$ are isomorphic to objects in
$\subcat$.
\end{Setup}

\begin{Remark}
\label{rmk:Krull-Schmidt}
By \cite[p.\ 52, thm.]{RingTame}, the conditions of Setup
\ref{set:blanket} imply that $\cat$ is a Krull-Schmidt category.  That
is, each object in $\cat$ is the coproduct of finitely many
indecomposable objects which are unique up to isomorphism.
\end{Remark}

\begin{Definition}
\label{def:arsubcat}
Let $\Sigma$ denote the suspension in $\cat$.  A distinguished
triangle 
\[
  \xymatrix{
  A \ar[r] & B \ar[r] & C \ar[r] & \Sigma A
           }
\]
with $A$, $B$, and $C$ in $\subcat$ is called an {\em Auslander-Reiten
triangle in $\subcat$} if it satisfies the following.
\begin{enumerate}

  \item  The morphism $\xymatrix{ C \ar[r] & \Sigma A }$ is non-zero.

\smallskip

  \item  If $A^{\prime}$ is in $\subcat$ then each morphism
         $\xymatrix{ A \ar[r] & A^{\prime} }$ which is not a section
         has a factorization 
\[
  \xymatrix{
    A \ar[r] \ar[d] & B \lefteqn{.} \ar[ld] \\
    A^{\prime} & \\
  }
\]

\smallskip

  \item  If $C^{\prime}$ is in $\subcat$ then each morphism
         $\xymatrix{ C^{\prime} \ar[r] & C }$ which is not a
         retraction has a factorization 
\[
  \xymatrix{
    & C^{\prime} \ar[d] \ar[ld] \\
    B \ar[r] & C \lefteqn{.} \\
  }
\]

\end{enumerate}
\end{Definition}

\begin{Definition}
\label{def:cover_and_envelope}
Let $X$ be an object of $\cat$.  A {\em $\subcat$-precover} of $X$ is
a morphism $\xymatrix{ A \ar[r] & X }$ with $A$ in $\subcat$, for
which each morphism $\xymatrix{ A^{\prime} \ar[r] & X }$ with
$A^{\prime}$ in $\subcat$ has a factorization
\[
  \xymatrix{
    & A \ar[d]  \\
    A^{\prime} \ar[r] \ar[ru] & X \lefteqn{.} \\
  }
\]

A {\em $\subcat$-cover} of $X$ is a $\subcat$-precover $\xymatrix{ A
\ar[r] & X }$ which is right-minimal, that is, satisfies that if
$\xymatrix{ A \ar[r] & A }$ is an endomorphism for which the
composition $\xymatrix{ A \ar[r] & A
\ar[r] & X }$ is equal to $\xymatrix{ A \ar[r] & X, }$ then
$\xymatrix{ A \ar[r] & A }$ is an isomorphism.

The notions of precover and cover have the dual notions of {\em
preenvelope} and {\em envelope}.
\end{Definition}

\section{Lemmas}
\label{sec:lemmas}

This section proves some lemmas which are needed later.  The following
is a triangulated version of Wakamatsu's Lemma; the proof is not too
far from the usual abelian case.

\begin{Lemma}
[Triangulated Wakamatsu's Lemma]
\label{lem:Wakamatsu}
Suppose that the morphism $\xymatrix{ A \ar[r]^{\alpha} & X }$ in
$\cat$ is a $\subcat$-cover and complete it to a distinguished
triangle $\xymatrix{ K \ar[r] & A \ar[r]^{\alpha} & X \ar[r] & . }$

Then $\Hom(C,\Sigma K) = 0$ for each $C$ in $\subcat$.
\end{Lemma}

\begin{proof}
There is a long exact sequence
\[
  \xymatrix{
    \Hom(C,A) \ar[rr]^{\Hom(C,\alpha)} & & \Hom(C,X) \ar[d] & & \\
    & & \Hom(C,\Sigma K) \ar[d] & & \\
    & & \Hom(C,\Sigma A) \ar[rr]_{\Hom(C,\Sigma(\alpha))} & & \Hom(C,\Sigma X).
           }
\]
The first map in the sequence is surjective because $\xymatrix{ A
\ar[r]^{\alpha} & X }$ is a $\subcat$-cover, so the second
map is zero.  I claim that the fourth map in the sequence is injective
whence the third map is zero.  This forces $\Hom(C,\Sigma K) = 0$ as
desired.

To see that the fourth map is injective, let $c$ in $\Hom(C,\Sigma A)$
have $\Hom(C,\Sigma(\alpha))(c) = \Sigma(\alpha) c = 0$.  This can be
interpreted as a commutative square
\[
  \xymatrix{
    C \ar[r]^c \ar@{=}[d] & \Sigma A \ar[d]^{\Sigma(\alpha)} \\
    C \ar[r]_0 & \Sigma X \lefteqn{,}
  }
\]
and completing the horizontal morphisms to distinguished triangles
gives a commutative diagram
\[
  \xymatrix{
    A \ar[r]^a \ar[d]_{\alpha} & B \ar[r] \ar[d]^{\beta} & C \ar[r]^c \ar@{=}[d] & \Sigma A \ar[d]^{\Sigma(\alpha)} \\
    X \ar[r]_x & M \ar[r] & C \ar[r]_0 & \Sigma X \lefteqn{,}
  }
\]
where $\beta$ exists by the axioms of triangulated categories.  In
particular, 
\[
  \beta a = x \alpha.
\]

The lower triangle is split so there is a morphism $\xymatrix{ M
\ar[r]^{m} & X }$ with
\[
  mx = \id_X.
\]
Consider the morphism $\xymatrix{ B \ar[r]^{m\beta} & X. }$  Since
$\subcat$ is closed under extensions, $B$ is in $\subcat$, and as
$\xymatrix{ A \ar[r]^{\alpha} & X }$ is a $\subcat$-cover, there is a
factorization
\[
  \xymatrix{
    A \ar[d]_{\alpha} & B \lefteqn{,} \ar[l]_b \ar[ld]^{m\beta} \\
    X & \\
  }
\]
that is,
\[
  \alpha b = m\beta.
\]
But now
\[
  \alpha \circ b a = m \beta a = mx\alpha = \id_X\alpha = \alpha,
\]
and as $\xymatrix{ A \ar[r]^{\alpha} & X }$ is a $\subcat$-cover, this
shows that $\xymatrix{ A \ar[r]^{ba} & A }$ is an automorphism.

Let $a^{\prime}$ be an inverse.  Then
\[
  \Sigma^{-1}(c)
  = \id_A \Sigma^{-1}(c)
  = a^{\prime} b a \Sigma^{-1}(c)
  = a^{\prime} b \circ 0
  = 0,
\]
where I have used that $a \Sigma^{-1}(c) = 0$ because $a$ and
$\Sigma^{-1}(c)$ are consecutive morphisms in a distinguished triangle.
The displayed equation shows $\Sigma^{-1}(c) = 0$, and hence $c = 0$ as
desired. 
\end{proof}

\begin{Lemma}
\label{lem:simple_socle}
Let $\xymatrix{ X \ar[r] & Y \ar[r]^{y} & Z \ar[r]^{d} & \Sigma X }$
be an Auslander-Reiten triangle in $\cat$.  View the abelian group
$\Hom(Z,\Sigma X)$ as a $\Hom(Z,Z)$-right-module via composition of
morphisms.

The socle of this module is simple and equal to the submodule
generated by $\xymatrix{ Z \ar[r]^{d} & \Sigma X. }$
\end{Lemma}

\begin{proof}
Let $M$ be a non-zero submodule of the $\Hom(Z,Z)$-right-module
$\Hom(Z,\Sigma X)$.  Pick a non-zero element $m$ in $M$, that is, $m$
is a non-zero morphism $\xymatrix{ Z \ar[r] & \Sigma X. }$  Then
$\xymatrix{ \Sigma^{-1}Z \ar[r]^-{\Sigma^{-1}(m)} & X }$ is also a
non-zero morphism, and it follows from \cite[lem.\ 3.3]{Krause2} that
there is a factorization
\[
  \xymatrix{
    & \Sigma^{-1}Z \ar[d]^{\Sigma^{-1}(m)} \\
    \Sigma^{-1}Z \ar[r]_-{\Sigma^{-1}(d)} \ar[ru]^{\Sigma^{-1}(z)} & X \lefteqn{,} \\
  }
\]
where the factoring morphism can clearly be taken to be of the form
$\Sigma^{-1}(z)$.  

Hence $\Sigma^{-1}(d) = \Sigma^{-1}(m) \Sigma^{-1}(z)$ and so $d =
mz$.  This means that in the $\Hom(Z,Z)$-right-module $\Hom(Z,\Sigma
X)$, the element $d$ is a multiple of $m$.  Hence $d$ is in $M$ and
consequently, the submodule of $\Hom(Z,\Sigma X)$ generated by $d$ is
contained in $M$.  Since $d$ is non-zero, so is the submodule
generated by $d$, and it follows that the socle of $\Hom(Z,\Sigma X)$
is equal to the submodule generated by $d$.

Now note that $\Hom(Z,Z)$ is a local ring by the dual of \cite[lem.\
2.3]{Krause1}.  Since the socle of $\Hom(Z,\Sigma X)$ is generated by
a single element, it will follow that it is simple if it is
annihilated by the Jacobson radical of $\Hom(Z,Z)$.  So let $r$ be in
the radical.  Then $r$ does not have a right-inverse by the dual of
\cite[prop.\ 15.15(e)]{AF}, and this is the same as to say that the
morphism $\xymatrix{ Z \ar[r]^{r} & Z }$ is not a retraction.  Hence
there is a factorization
\[
  \xymatrix{
    & Z \ar[d]^r \ar[dl]_{r^{\prime}} \\
    Y \ar[r]_y & Z \\
  }
\]
with $r = yr^{\prime}$, and so
\[
  dr = dyr^{\prime} = 0 \circ r^{\prime} = 0
\]
as desired, because $d$ and $y$ are consecutive morphisms in a
distinguished triangle whence $dy = 0$.
\end{proof}

\begin{Lemma}
\label{lem:indecomposable}
Let $C$ be in $\subcat$ and let $\xymatrix{ X \ar[r] & Y \ar[r] & C
\ar[r]^{d} & \Sigma X }$ be an Auslander-Reiten triangle in
$\cat$.

Suppose that $\xymatrix{ A \ar[r]^{\alpha} & X }$ is a
$\subcat$-cover.  Then $A$ is either zero or indecomposable.
\end{Lemma}

\begin{proof}
Suppose that $A$ is not zero.  Recall from Remark
\ref{rmk:Krull-Schmidt} that $\cat$ is a Krull-Schmidt category.  Let
$A_i$ be an indecomposable direct summand of $A$, let $\xymatrix{ A_i
\ar@{^{(}->}[r] & A }$ be the inclusion of $A_i$ into $A$, and denote
by $\alpha_i$ the composition $\xymatrix{ A_i \ar@{^{(}->}[r] & A
\ar[r]^{\alpha} & X. }$  Since $\xymatrix{ A \ar[r]^{\alpha} & X }$ is
a cover, it is clear that $\alpha_i \not= 0$.  It follows from
\cite[lem.\ 3.3]{Krause2} that there is a factorization
\[
  \xymatrix{
    & A_i \ar[d]^{\alpha_i} \\
    \Sigma^{-1}C \ar[r]_-{\Sigma^{-1}(d)} \ar[ur]^s & X \lefteqn{.} \\
  }
\]
Since $d \not= 0$, it follows that $\Sigma^{-1}(d) \not= 0$ and hence $s
\not= 0$.

So each indecomposable direct summand $A_i$ of $A$ permits a non-zero
morphism $\xymatrix{ \Sigma^{-1}C \ar[r] & A_i. }$  To finish the
proof, let me show that at most one indecomposable summand can permit
such a morphism:

The $\subcat$-cover $\xymatrix{ A \ar[r]^{\alpha} & X }$ can be
completed to give a distinguished triangle $\xymatrix{ K \ar[r] & A
\ar[r]^{\alpha} & X \ar[r] & , }$ which can be turned to a
distinguished triangle $\xymatrix{ \Sigma K \ar[r] & \Sigma A \ar[r] &
\Sigma X \ar[r] & . }$  I have $\Hom(C,\Sigma K) = 0$ by Lemma
\ref{lem:Wakamatsu}, and hence the homomorphism
\[
  \xymatrix{
  \Hom(C,\Sigma A) \ar[r] & \Hom(C,\Sigma X)
           }
\]
is injective.  Viewing this as a homomorphism of finite dimensional
right-modules over the finite dimensional $k$-algebra $\Hom(C,C)$, the
target has simple socle by Lemma \ref{lem:simple_socle}.  Hence the
image is either zero or indecomposable, and since the homomorphism is
injective, the same is true for the source $\Hom(C,\Sigma A)$.  So if
$A$ splits as $A \cong A_1 \amalg \cdots \amalg A_n$, there can be at
most one $i$ for which $\Hom(C,\Sigma A_i) \cong
\Hom(\Sigma^{-1}C,A_i)$ is non-zero.

That is, there is at most one indecomposable summand $A_i$ of $A$
which permits a non-zero morphism $\xymatrix{ \Sigma^{-1}C \ar[r] &
A_i, }$ as desired.
\end{proof}

%\begin{Lemma}
%\label{lem:precover_and_cover}
%If $A \rightarrow X$ is a $\subcat$-precover, then $A$ has a direct
%summand $A_1$ such that $A_1 \rightarrow X$ is a $\subcat$-cover.
%\end{Lemma}

The easy proof of the following lemma is left to the reader.  As a
hint, local rings are precisely the rings where the set of elements
without a left-inverse is closed under addition, and also precisely
the rings where the set of elements without a right-inverse is closed
under addition, cf.\ \cite[prop.\ 15.15]{AF}.  See also \cite[lem.\
2.3]{Krause1}.

\begin{Lemma}
\label{lem:local_End}
Let $\xymatrix{ A \ar[r] & B \ar[r] & C \ar[r] & {} }$ be an
Auslander-Reiten triangle in the subcategory $\subcat$.

Then $A$ and $C$ have local endomorphism rings.  In particular, $A$
and $C$ are indecomposable.
\end{Lemma}

%THE FOLLOWING is a complete proof.
%\begin{proof}
%It is enough to show that the endomorphism ring of $A$ is local, the
%case of $C$ being dual.
%
%By \cite[prop.\ 15.15]{AF}, I must show that the set of elements of
%the endomorphism ring without left-inverses is closed under addition.
%These elements are precisely the morphisms $A \rightarrow A$ which are
%not sections.
%
%Given two such elements, $a_1$ and $a_2$, they both factor through $A
%\stackrel{a}{\rightarrow} B$, and hence $s = a_1 + a_2$ also factors
%through $A \stackrel{a}{\rightarrow} B$, that is, $s = ba$ for
%some $B \stackrel{b}{\rightarrow} A$.
%
%But then the sum $s$ cannot have a left-inverse $a^{\prime}$, and this
%proves the Lemma.  For suppose that $a^{\prime}$ existed.  Then
%$a^{\prime}s = \id_A$, and then the morphism 
%$\xymatrix{
%\Sigma^{-1}C \ar[r]^-{\Sigma^{-1}(\frakd)} & A \\
%}$
%would satisfy
%\[
%  \Sigma^{-1}(\frakd) 
%    = \id_A \Sigma^{-1}(\frakd)
%    = a^{\prime} s \Sigma^{-1}(\frakd)
%    = a^{\prime} b a \Sigma^{-1}(\frakd)
%    = a^{\prime} b \circ 0
%    = 0,
%\]
%where I have used that $a \Sigma^{-1}(\frakd) = 0$ because $a$ and
%$\Sigma^{-1}(\frakd)$ are consecutive morphisms in a distinguished
%triangle.  But this contradicts $\frakd \not= 0$.
%\end{proof}

Finally, an elementary fact of linear algebra.

\begin{Lemma}
\label{lem:bilinear}
Let $U$ and $V$ be finite dimensional $k$-vector spaces and let
\[
  \xymatrix{
  q : U \times V \ar[r] & k
           }
\]
be a bilinear map such that for each non-zero $u$ in $U$, there
exists a $v$ in $V$ with $q(u,v) \not= 0$.

Then for each linear map $\xymatrix{ \varphi : U \ar[r] & k, }$ there
exists a $v$ in $V$ such that
\[
  \varphi(-) = q(-,v).
\]
\end{Lemma}

\section{Main Theorem}
\label{sec:main}

The following is the main theorem of this paper.  Recall the
triangulated category $\cat$ and the subcategory $\subcat$ from Setup
\ref{set:blanket}.

\begin{Theorem}
\label{thm:main}
Let $C$ be in $\subcat$ and suppose that there exists an $A^{\prime}$
in $\subcat$ and a non-zero morphism $\xymatrix{ C \ar[r]^{c} &
\Sigma A^{\prime}. }$

Let
\begin{equation}
\label{equ:d}
  \xymatrix{
    X \ar[r] & Y \ar[r]^{y} & C \ar[r]^{d} & \Sigma X
           }
\end{equation}
be an Auslander-Reiten triangle in $\cat$.  Then the following are
equivalent.
\begin{enumerate}

  \item  $X$ has a $\subcat$-cover of the form $\xymatrix{ A
         \ar[r]^{\alpha} & X. }$

\smallskip

  \item  There is an Auslander-Reiten triangle in $\subcat$,
\[
  \xymatrix{
  A \ar[r] & B \ar[r] & C \ar[r] & .
           }
\]

\end{enumerate}
\end{Theorem}

Before the proof, let me remark that the existence of a non-zero
morphism $\xymatrix{ C \ar[r] & \Sigma A^{\prime} }$ is a necessary
condition for the existence of an Auslander-Reiten triangle
$\xymatrix{ A \ar[r] & B \ar[r] & C \ar[r] & {} }$ in $\subcat$, for
if the triangle exists then its connecting morphism $\xymatrix{ C
\ar[r] & \Sigma A }$ is non-zero.

\begin{proof}
(i) $\Rightarrow$ (ii).  The morphism $\xymatrix{ C \ar[r]^{c}
& \Sigma A^{\prime} }$ is non-zero, so com\-ple\-ting it to a
distinguished triangle
\[
  \xymatrix{
  A^{\prime} \ar[r] & B^{\prime} \ar[r]^{b^{\prime}} & C 
    \ar[r]^{c} & \Sigma A^{\prime}
           }
\]
gives that $b^{\prime}$ is not a retraction.  Hence the definition of
Auslander-Reiten triangles means that $\beta^{\prime}$ exists in the
following commutative diagram,
\[
  \xymatrix
  {
    A^{\prime} \ar[r] \ar[d]_{\alpha^{\prime}_1} & B^{\prime} \ar[r]^{b^{\prime}} \ar[d]_{\beta^{\prime}} & C \ar[r]^{c} \ar@{=}[d] & \Sigma A^{\prime} \ar[d]^{\Sigma(\alpha^{\prime}_1)} \\
    X \ar[r] & Y \ar[r]_y & C \ar[r]_{d} & \Sigma X \lefteqn{,} \\
  }
\]
and $\alpha^{\prime}_1$ exists by the axioms of triangulated
categories.  Note that
\[
  \Sigma(\alpha_1^{\prime})c = d.
\]
Since $\xymatrix{ A \ar[r]^{\alpha} & X }$ is a $\subcat$-cover, there
is a factorisation
\[
  \xymatrix{
    & A^{\prime} \ar[d]^{\alpha_1^{\prime}} \ar[dl]_{\alpha_2^{\prime}} \\
    A \ar[r]_{\alpha} & X \\
  }
\]
with
\[
  \alpha\alpha_2^{\prime} = \alpha_1^{\prime}.
\]
Consider the morphism $\xymatrix{ C
\ar[r]^{\Sigma(\alpha^{\prime}_2)c} & \Sigma A }$ and complete it to a
distinguished triangle
\begin{equation}
\label{equ:a}
  \xymatrix
  {
    A \ar[r] & B \ar[r]^b & C \ar[r]^-{\Sigma(\alpha^{\prime}_2)c} & \Sigma A \lefteqn{.}
  }
\end{equation}
I claim that this is an Auslander-Reiten triangle in $\subcat$.

To see so, note first that since $A$ and $C$ are in $\subcat$, so is
$B$, because $\subcat$ is closed under extensions.  Let me next verify
that the conditions of Definition \ref{def:arsubcat} hold for the
distinguished triangle \eqref{equ:a}.

For condition (i), note that
\begin{equation}
\label{equ:c}
  \Sigma(\alpha) \circ \Sigma(\alpha^{\prime}_2)c 
  = \Sigma(\alpha\alpha^{\prime}_2)c 
  = \Sigma(\alpha^{\prime}_1)c 
  = d,
\end{equation}
and since the connecting morphism $d$ of the Auslander-Reiten triangle
\eqref{equ:d} is non-zero, it follows that $\Sigma(\alpha^{\prime}_2)c
\not= 0$.

Observe for later use that in particular, the target $\Sigma A$ of
$\Sigma(\alpha^{\prime}_2)c$ must be non-zero, whence also
\begin{equation}
\label{equ:e}
  A \not\cong 0.
\end{equation}

For condition (ii), consider the commutative square
\[
  \xymatrix
  {
    C \ar[r]^-{\Sigma(\alpha^{\prime}_2)c} \ar@{=}[d] & \Sigma A \ar[d]^{\Sigma(\alpha)} \\
    C \ar[r]_{d} & \Sigma X \\
  }
\]
which exists by equation \eqref{equ:c}.  Using the octahedral axiom,
this can be extended to a commutative diagram
\[
  \xymatrix
  {
    A \ar[r] \ar[d]_{\alpha} & B \ar[r]^b \ar[d]_{\beta} & C \ar[r]^-{\Sigma(\alpha^{\prime}_2)c} \ar@{=}[d] & \Sigma A \ar[d]^{\Sigma(\alpha)} \\
    X \ar[r] \ar[d] & Y \ar[r]_y \ar[d]_{\upsilon} & C \ar[r]_{d} \ar[d] & \Sigma X \ar[d] \\
    \Sigma K \ar@{=}[r] \ar[d] & \Sigma K \ar[r] \ar[d] & 0 \ar[r] \ar[d] & \Sigma^2 K \ar[d] \\
    \Sigma A \ar[r] & \Sigma B \ar[r] & \Sigma C \ar[r] & \Sigma^2 A \\
  }
\]
where the first row is the distinguished triangle \eqref{equ:a}, the
second row is the Auslander-Reiten triangle \eqref{equ:d}, and each
row and each column is a distinguished triangle.  

Now, given $C^{\prime}$ in $\subcat$, each morphism $\xymatrix{
C^{\prime} \ar[r] & C }$ which is not a retraction has a factorization
\[
  \xymatrix{
    & C^{\prime} \ar[d] \ar[ld] \\
    Y \ar[r]_y & C \\
  }
\]
because the second row in the diagram is an Auslander-Reiten triangle.
The composition $\xymatrix{ C^{\prime} \ar[r] & Y \ar[r]^{\upsilon} &
\Sigma K }$ is zero by Lemma \ref{lem:Wakamatsu}, so the morphism
$\xymatrix{ C^{\prime} \ar[r] & Y }$ again has a factorization
\[
  \xymatrix{
    B \ar[d]_{\beta} & C^{\prime} \lefteqn{.} \ar[l] \ar[ld] \\
    Y & \\
  }
\]
So the morphism $\xymatrix { C^{\prime} \ar[r] & C }$ has been
factored as
\[
  \xymatrix{
    C^{\prime} \ar[r] & B \ar[r]^{\beta} & Y \ar[r]^{y} & C,
           }
\]
that is, as $\xymatrix{ C^{\prime} \ar[r] & B \ar[r]^{b} & C }$,
verifying condition (ii) of Definition \ref{def:arsubcat}. 

To see that condition (iii) of Definition \ref{def:arsubcat} holds for
the distinguished triangle \eqref{equ:a}, note that $A$ is non-zero by
equation \eqref{equ:e}, hence indecomposable by Lemma
\ref{lem:indecomposable}, and so $\Hom(A,A)$ is a local ring because
$\cat$ has finite dimensional $\Hom$ spaces and split idempotents.

But then $\Hom(\Sigma A,\Sigma A)$ is also a local ring since it is
isomorphic to $\Hom(A,A)$.  Consequently, \cite[lem.\
2.4]{Krause1} implies that in \eqref{equ:a},
the third morphism $\xymatrix{ C
\ar[r]^{\Sigma(\alpha^{\prime}_2)c} & \Sigma A }$
is left minimal, and then \cite[lem.\ 2.5]{Krause1}
implies that the second morphism $\xymatrix{ B \ar[r]^{b} & C }$ of
\eqref{equ:a} is right minimal.  Since I have already proved that the
distinguished triangle \eqref{equ:a} satisfies condition (ii) of
Definition \ref{def:arsubcat}, condition (iii) can now be proved by
the same method as in \cite[proof of lem.\ 2.6,
(2)$\Rightarrow$(3)]{Krause1}.

\smallskip

(ii) $\Rightarrow$ (i).  Let
\[
  \xymatrix{
    A \ar[r] & B \ar[r]^{b} & C \ar[r]^{\frakd} & \Sigma A
           }
\]
be an Auslander-Reiten triangle in $\subcat$.  Since $\frakd$ is
non-zero, $\xymatrix{ B \ar[r]^{b} & C }$ is not a retraction, and so 
by the definition of Auslander-Reiten triangles and the axioms of
triangulated categories there is a commutative diagram
\[
  \xymatrix
  {
    \Sigma^{-1}C \ar[r]^-{\Sigma^{-1}(\frakd)} \ar@{=}[d] & A \ar[r] \ar[d]^{\alpha} & B \ar[r]^b \ar[d] & C \ar@{=}[d] \\
    \Sigma^{-1}C \ar[r]_-{\Sigma^{-1}(d)} & X \ar[r] & Y \ar[r]_y & C \lefteqn{;} \\
  }
\]
in particular,
\[
  \alpha\Sigma^{-1}(\frakd) = \Sigma^{-1}(d).
\]
I claim that $\xymatrix{ A \ar[r]^{\alpha} & X }$ is a
$\subcat$-cover.  In fact, $A$ has local endomorphism ring by Lemma
\ref{lem:local_End}, so $\xymatrix{ A \ar[r]^{\alpha} & X }$ is right
minimal by the dual of \cite[lem.\ 2.4]{Krause1}.  Hence it is enough
to show that $\xymatrix{ A \ar[r]^{\alpha} & X }$ is a
$\subcat$-precover, cf.\ Definition \ref{def:cover_and_envelope}.

To prove this, pick a linear map
\[
  \xymatrix{
    \psi : \Hom(\Sigma^{-1}C,X) \ar[r] & k
           }
\]
which satisfies $\psi(\Sigma^{-1}(d)) \not= 0$.  Given a non-zero
$A^{\prime}$ in $\subcat$, define a bilinear map
\begin{eqnarray*}
  & \xymatrix{
      q : \Hom(\Sigma^{-1}C,A^{\prime}) \times \Hom(A^{\prime},A)
      \ar[r] & k,
             } & \\
  & \xymatrix{ q(s,a^{\prime}) = \psi(\alpha a^{\prime}s). } &
\end{eqnarray*}

Let me first show that Lemma \ref{lem:bilinear} applies to $q$.  Let
$s$ be a non-zero element in $\Hom(\Sigma^{-1}C,A^{\prime})$.
Construct a distinguished triangle
\[
  \xymatrix{
  A^{\prime} \ar[r] & B^{\prime} \ar[r]^{b^{\prime}} & C
    \ar[r]^{\Sigma(s)} & \Sigma A^{\prime}.
           }
\]
Since $\subcat$ is closed under extensions, $B^{\prime}$ is in
$\subcat$.  Since $s$ is non-zero, so is $\Sigma(s)$ and hence
$\xymatrix{ B^{\prime} \ar[r]^{b^{\prime}} & C }$ is not a
retraction.  So there is a factorization
\[
  \xymatrix
  {
    & B^{\prime} \ar[ld]_{\beta^{\prime}} \ar[d]^{b^{\prime}} \\
    B \ar[r]_b & C \lefteqn{.} \\
  }
\]
Hence
\[
  \frakd b^{\prime} = \frakd b \beta^{\prime} = 0 \circ \beta^{\prime} = 0,
\]
where I have used that $\frakd b = 0$ since $\frakd$ and $b$ are
consecutive morphisms in a distinguished triangle.  So there is a
factorization
\[
  \xymatrix{
    B^{\prime} \ar[r]^{b^{\prime}} & C \ar[r]^{\Sigma(s)} \ar[d]_{\frakd} & \Sigma A^{\prime} \ar[ld]^{\Sigma(a^{\prime})}\lefteqn{,} \\
    & \Sigma A & \\
  }
\]
where the factoring morphism can clearly be taken to be of the form
$\Sigma(a^{\prime})$.  Since $\Sigma(a^{\prime}) \Sigma(s) = \frakd$ I
get
\[
  a^{\prime}s = \Sigma^{-1}(\frakd)
\]
and so the element $a^{\prime}$ in $\Hom(A^{\prime},A)$ satisfies
\[
  q(s,a^{\prime}) = \psi(\alpha a^{\prime}s) 
                   = \psi(\alpha \Sigma^{-1}(\frakd))
                   = \psi(\Sigma^{-1}(d)) \not= 0,
\]
proving that Lemma \ref{lem:bilinear} applies to $q$ as desired.

Now let $\xymatrix{ A^{\prime} \ar[r]^{\alpha^{\prime}} & X }$ be
given.  I want to factor it through $\xymatrix{ A \ar[r]^{\alpha} & X
}$ in order to prove that $\xymatrix{ A \ar[r]^{\alpha} & X }$ is a
$\subcat$-precover.  The case $A^{\prime} = 0$ is trivial, so suppose
$A^{\prime} \not= 0$.  Consider the linear map
\begin{eqnarray*}
  & \xymatrix{
      \varphi : \Hom(\Sigma^{-1}C,A^{\prime}) \ar[r] & k,
             } & \\
  & \xymatrix{ \varphi(s) = \psi(\alpha^{\prime}s). } &
\end{eqnarray*}
By Lemma \ref{lem:bilinear}, there exists an element $a^{\prime}$ in
$\Hom(A^{\prime},A)$ such that $\varphi(-) = q(-,a^{\prime})$, and by
the definition of $\varphi$ and $q$, this says
\begin{equation}
\label{equ:b}
  \psi(\alpha^{\prime}s) = \psi(\alpha a^{\prime}s) \;\;
  \mbox{for each $s$.}
\end{equation}

However, \cite[proof of prop.\ I.2.3]{RVdB} shows that the bilinear
map
\begin{eqnarray*}
  & \xymatrix{
      \Hom(\Sigma^{-1}C,A^{\prime}) \times \Hom(A^{\prime},X) \ar[r] & k,
             } & \\
  & \xymatrix{ (s,\fa) \mapsto \psi(\fa s) } &
\end{eqnarray*}
is non-degenerate, so equation \eqref{equ:b} implies
\[
  \alpha^{\prime} = \alpha a^{\prime}.
\]
That is, I have obtained a factorization
\[
  \xymatrix{
    & A \ar[d]^{\alpha} \\
    A^{\prime} \ar[r]_{\alpha^{\prime}} \ar[ru]^{a^{\prime}} & X \\
  }
\]
as desired.
\end{proof}

The dual result can be proved dually.  Let me state it for easy
reference.

\begin{Theorem}
Let $A$ be in $\subcat$ and suppose that there exists a $C^{\prime}$ in
$\subcat$ and a non-zero morphism $\xymatrix{ C^{\prime} \ar[r] &
\Sigma A. }$

Let
\[
  \xymatrix{
    A \ar[r] & Y \ar[r] & Z \ar[r] & {}
           }
\]
be an Auslander-Reiten triangle in $\cat$.  Then the following are
equivalent.
\begin{enumerate}

  \item  $Z$ has a $\subcat$-envelope of the form $\xymatrix{ Z \ar[r]
         & C. }$

\smallskip

  \item  There is an Auslander-Reiten triangle in $\subcat$,
\[
  \xymatrix{
    A \ar[r] & B \ar[r] & C \ar[r] & .
           }
\]

\end{enumerate}
\end{Theorem}

\section{Existence of Auslander-Reiten sequences}
\label{sec:example}

This section uses the results of Section \ref{sec:main} to give a new
proof of the existence of Aus\-lan\-der-Reiten sequences over a finite
dimensional $k$-algebra $\Lambda$.  The idea is to consider the module
category $\mod(\Lambda)$ as a subcategory of a suitable triangulated
category of complexes of injective $\Lambda$-modules.

I thank Henning Krause for showing me how to use the methods from
\cite{Krause3} and \cite{KrauseLe} to remove the assumption
$\gldim \Lambda < \infty$.  One of the important tools is the
Auslander-Reiten triangles in the homotopy category of complexes of
injective $\Lambda$-modules described in \cite{KrauseLe}.  Note that
\cite{KrauseLe} already uses these to get Auslander-Reiten sequences;
the present section offers an alternative way to do so via the results
of Section \ref{sec:main}.

\begin{Lemma}
\label{lem:covers}
Let $\xymatrix{ A \ar[r]^{\alpha} & X }$ be a $\subcat$-precover.  
Then $A$ has a direct summand $A_1$ with inclusion $\xymatrix{ A_1
\ar@{^{(}->}[r] & A }$ such that the composition
\[
  \xymatrix{ A_1 \ar@{^{(}->}[r] & A \ar[r]^{\alpha} & X }
\]
is a $\subcat$-cover.
\end{Lemma}

\begin{proof}
Consider the category $\finpres(\cat)$ of finitely presented
contravariant additive functors from $\cat$ to the category of abelian
groups.  This is also known as the Freyd category; see
\cite[chp.\ 5]{Neemanbook} for a survey of the theory.  

The category $\finpres(\cat)$ is abelian, and the additive functor
\[
  \xymatrix{ \cat \ar[r] & \finpres(\cat), } \;\;\;
  \xymatrix{ X \ar@{|->}[r] & \Hom(-,X) }
\] 
permits me to view $\cat$ as the full subcategory of projective
objects of $\finpres(\cat)$.

To get $A_1$, view the $\subcat$-precover $\xymatrix{ A
\ar[r]^{\alpha} & X }$ as a morphism in $\finpres(\cat)$ and factor it
into an epimorphism followed by a monomorphism,
\[
  \xymatrix{
    A \ar[rr]^{\alpha} \ar@{->>}[rd] & & X \lefteqn{.} \\
    & \Image(\alpha) \ar@{^{(}->}[ru] & \\
  }
\]
It is not hard to verify that the direct summand $A_1$ can be obtained
as a projective cover $\xymatrix{ A_1 \ar@{->>}[r] &
\Image(\alpha) }$ in $\finpres(\cat)$ which exists by \cite[thm.\
4.12]{AusArtI}. 
\end{proof}

\begin{Setup}
\label{set:modules}
Let $\Lambda$ be a finite dimensional $k$-algebra over the field $k$.

Consider $\sK(\Inj\:\Lambda)$, the homotopy category of complexes of
injective $\Lambda$-left-modules, and $\dual\!\Lambda =
\Hom_k(\Lambda,k)$, the $k$-linear dual of $\Lambda$ which is a
$\Lambda$-bi-module.

Let $\cat$ be the full subcategory of $\sK(\Inj\:\Lambda)$ consisting
of complexes $X$ for which each $X^i$ is finitely generated, and where
$\H^i(X) = 0$ for $i \gg 0$ and
$\H^i(\Hom_{\Lambda}(\dual\!\Lambda,X)) = 0$ for $i \ll 0$.

Let $\subcat$ be the full subcategory of $\cat$ which consists of
injective resolutions of finitely generated $\Lambda$-left-modules.
That is, $\subcat$ consists of the complexes in $\cat$ of the form
\[
  \xymatrix{
    \cdots \ar[r] & 0 \ar[r] & A^0 \ar[r] & A^1 \ar[r] & \cdots
           }
\]
where all other cohomology groups than $\H^0(A)$ are zero.
\end{Setup}

\begin{Lemma}
The categories $\cat$ and $\subcat$ of Setup \ref{set:modules} satisfy
the requirements of Setup \ref{set:blanket}.
\end{Lemma}

\begin{proof}
It is easy to see that $\cat$ is a skeletally small $k$-linear
triangulated category.

To see that $\cat$ has split idempotents, note first that
$\sK(\Inj\:\Lambda)$ does.  Namely, by \cite[prop.\
1.6.8]{Neemanbook}, it is enough to see that $\sK(\Inj\:\Lambda)$ has
set indexed coproducts, and this is clear because the coproduct of a
set indexed family of injective $\Lambda$-left-modules is again
injective.

Now let $X$ be in $\cat$ and let $\xymatrix{X \ar[r]^e & X}$ be an
idempotent.  Then there is a splitting of $e$ in $\sK(\Inj\:\Lambda)$,
that is, a diagram
\begin{equation}
\label{equ:aa}
  \xymatrix{
    X_1 \ar@<1ex>[r]^{x_1} & X \ar@<1ex>[l]^x
           }
\end{equation}
in $\sK(\Inj\:\Lambda)$ with $x_1x = e$ and $xx_1 = \id_{X_1}$.  By
\cite[App.\ B]{Krause3}, I can assume that $X_1$ is a so-called
homotopically minimal complex, and then $xx_1 = \id_{X_1}$ implies
that the chain map underlying the homotopy class $xx_1$ is invertible.
That is, each component of the chain map is bijective, so $X_1^i$ is
isomorphic to a direct summand of $X^i$ for each $i$.  But each $X^i$
is finitely generated by the definition of $\cat$, so the same must be
true for each $X_1^i$.

Since it is clear that the direct summand $X_1$ also inherits the
properties $\H^i(X_1) = 0$ for $i \gg 0$ and
$\H^i(\Hom_{\Lambda}(\dual\!\Lambda,X_1)) = 0$ for $i \ll 0$ from $X$,
it follows that $X_1$ is in $\cat$, and so \eqref{equ:aa} is in fact a
splitting of $e$ in $\cat$.  Hence $\cat$ has split idempotents.

To see that $\cat$ has finite dimensional $\Hom$ spaces over $k$, note
that a complex $X$ is the mapping cone of the chain map
\[
  \xymatrix{
    \cdots \ar[r] & X^{-2} \ar[r] \ar[d] & X^{-1} \ar[r] \ar[d]^{\partial^{-1}} & 0 \ar[r] \ar[d] & \cdots \\
    \cdots \ar[r] & 0 \ar[r] & X^0 \ar[r] & X^1 \ar[r] & \cdots \lefteqn{,} \\
  }
\]
where $\partial$ is the differential of $X$.  Using the notation
$X^{\sqsupset}$ and $X^{\sqsubset}$ for the two displayed complexes,
this implies that there is a distinguished triangle
$\xymatrix{X^{\sqsupset} \ar[r] & X^{\sqsubset} \ar[r] & X \ar[r] &
{}}$ which can be turned to a distinguished triangle
\[
  \xymatrix{
    X^{\sqsubset} \ar[r] & X \ar[r] & \Sigma X^{\sqsupset} \ar[r] & .
           }
\]
If $X$ is in $\cat$ then it is easy to check that the complexes
$X^{\sqsubset}$ and $\Sigma X^{\sqsupset}$ are also in $\cat$.  So
using an obvious notion of left- and right-boundedness for complexes,
the last triangle shows that each $X$ in $\cat$ is the extension
of a left-bounded complex $X^{\sqsubset}$ and a right-bounded complex
$\Sigma X^{\sqsupset}$, both of which are in $\cat$.

Hence it is sufficient to show that $\Hom$ spaces between left- and
right-bounded complexes in $\cat$ are finite dimensional.  This is
clear for the $\Hom$ space from a left-bounded to a right-bounded
complex, since the non-zero parts of the complexes only overlap in
finitely many degrees, and since the modules in the complexes are
finitely generated over $\Lambda$, hence finite dimensional over $k$.
The same argument applies to the $\Hom$ space from a right-bounded to
a left-bounded complex.

If $X$ and $Y$ are both left-bounded complexes in $\cat$, then $X$ and
$Y$ are left-bounded complexes consisting of injective modules, that
is, they are injective resolutions.  By the definition of $\cat$, the
cohomology of $X$ and $Y$ is right-bounded and finitely generated over
the finite dimensional $k$-algebra £$\Lambda$.  This implies that the
$\Hom$ space $\Hom_{\Dsmall(\Lambda)}(X,Y)$ in the derived category
$\D(\Lambda)$ is finite dimensional, and so the same is true for
$\Hom_{\catsmall}(X,Y)$ since
\[
  \Hom_{\Dsmall(\Lambda)}(X,Y)
    \cong \Hom_{\mbox{\sss K}(\mbox{\sss Inj}\:\Lambda)}(X,Y)
    \cong \Hom_{\catsmall}(X,Y),
\]
where the first $\cong$ is because $Y$ is an injective resolution and
the second $\cong$ is because $\cat$ is a full subcategory of
$\sK(\Inj\:\Lambda)$.

If $X$ and $Y$ are both right-bounded complexes in $\cat$, then
\[
  P = \Hom_{\Lambda}(\dual\!\Lambda,X) \;\; \mbox{and} \;\;
  Q = \Hom_{\Lambda}(\dual\!\Lambda,Y)
\]
are right-bounded complexes consisting of projective
$\Lambda$-left-modules, that is, $P$ and $Q$ are projective
resolutions.  By the definition of $\cat$, the cohomology of $P$ and
$Q$ is left-bounded and finitely generated.  This implies that
$\Hom_{\Dsmall(\Lambda)}(P,Q)$ is finite dimensional, and so the same
is true for $\Hom_{\catsmall}(X,Y)$ since
\begin{align*}
  \Hom_{\Dsmall(\Lambda)}(P,Q)
  & \stackrel{\rm (a)}{\cong} \Hom_{\mbox{\sss K}(\Lambda)}(P,Q) \\
  & = \Hom_{\mbox{\sss K}(\Lambda)}
          (\Hom_{\Lambda}(\dual\!\Lambda,X),\Hom_{\Lambda}(\dual\!\Lambda,Y))\\
  & \stackrel{\rm (b)}{\cong} \Hom_{\mbox{\sss K}(\mbox{\sss Inj}\:\Lambda)}(X,Y) \\
  & \stackrel{\rm (c)}{\cong} \Hom_{\catsmall}(X,Y).
\end{align*}
Here $\sK(\Lambda)$ denotes the homotopy category of complexes of
$\Lambda$-left-mo\-du\-les, (a) is because $P$ is a projective
resolution, (b) is because the functor
$\Hom_{\Lambda}(\dual\!\Lambda,-))$ is an equivalence from the
injective $\Lambda$-left-modules to the projective
$\Lambda$-left-modules, and (c) is because $\cat$ is a full
subcategory of $\sK(\Inj\:\Lambda)$.  This shows that $\cat$ satisfies
the requirements of Setup \ref{set:blanket}.

Finally, it is easy to see that $\subcat$ is closed under extensions
and direct summands.
\end{proof}

\begin{Remark}
\label{rmk:mod}
$\subcat$ is the full subcategory of $\cat$ which consists of
injective resolutions of finitely generated $\Lambda$-left-modules.
Since the morphisms in $\subcat$ are homotopy classes of chain maps,
it is classical that $\subcat$ is e\-qui\-va\-lent to $\mod(\Lambda)$,
the category of finitely generated $\Lambda$-left-modules.  An
equivalence can be constructed by letting a module correspond to one
of its injective resolutions.

Moreover, let
\[
  \xymatrix{
    A \ar[r] & B \ar[r] & C \ar[r] & {}
           }
\]
be an Auslander-Reiten triangle in the subcategory $\subcat$ of
$\cat$, in the sense of Definition \ref{def:arsubcat}.  Then
the cohomology long exact sequence gives a short exact sequence
\[
  \xymatrix{
  0 \ar[r] & \H^0(A) \ar[r] & \H^0(B) \ar[r] & \H^0(C) \ar[r] & 0,
           }
\]
and it is easy to see that this is an Auslander-Reiten sequence in the
abelian category $\mod(\Lambda)$.
\end{Remark}

\begin{Lemma}
\label{lem:mod_covering}
$\subcat$ is a covering class in $\cat$.
\end{Lemma}

\begin{proof} 
By Lemma \ref{lem:covers} it is enough to show that $\subcat$ is
precovering in $\cat$.  But given $X$ in $\cat$, let $\partial$ be the
differential of the complex $X$ and let $A$ be an injective resolution
of $\Ker \partial^0$.  It is then easy to check that the canonical
chain map
\[
  \xymatrix{
    \cdots \ar[r] & 0 \ar[r] \ar[d] & A^0 \ar[r] \ar[d] & A^1 \ar[r] \ar[d] & \cdots \\
    \cdots \ar[r] & X^{-1} \ar[r] & X^0 \ar[r] \ar[r] & X^1 \ar[r] & \cdots \\
  }
\]
gives a $\subcat$-precover of $X$ in $\cat$.
\end{proof}

\begin{Theorem}
Let $M$ and $U$ be indecomposable finitely generated
$\Lambda$-left-modules and suppose that $M$ is not projective and that
$U$ is not injective.

Then there are Auslander-Reiten sequences in $\mod(\Lambda)$,
\[
  \xymatrix{
    0 \ar[r] & K \ar[r] & L \ar[r] & M \ar[r] & 0
           }
\]
and
\[
  \xymatrix{
  0 \ar[r] & U \ar[r] & V \ar[r] & W \ar[r] & 0.
           }
\]
\end{Theorem}

\begin{proof}
It is enough to show that the first of these sequences exists, because
the second one can be obtained using the dualization functor $\dual(-)
= \Hom_k(-,k)$.

Since $M$ is not a projective module, there exists a $K^{\prime}$ in
$\mod(\Lambda)$ such that $\Ext_{\Lambda}^1(M,K^{\prime}) \not= 0$.
Let $C$ and $A^{\prime}$ be injective resolutions of $M$ and
$K^{\prime}$.  Then $C$ and $A^{\prime}$ are in $\subcat$, and since
\[
  \Ext_{\Lambda}^1(M,K^{\prime}) 
  \cong \Hom_{\mbox{\sss K}(\mbox{\sss Inj}\:\Lambda)}(C,\Sigma A^{\prime})
  \cong \Hom_{\catsmall}(C,\Sigma A^{\prime}),
\]
it follows that there is a non-zero morphism $\xymatrix{ C \ar[r] &
\Sigma A^{\prime} }$ in $\cat$.

The complex $C$ is a compact object of $\sK(\Inj\:\Lambda)$ by
\cite[lem.\ 2.1]{Krause3}, and since $M$ is an indecomposable module,
$C$ is an indecomposable object of $\sK(\Inj\:\Lambda)$.  Hence by
\cite[thm.\ 6.3]{KrauseLe}, there is an Auslander-Reiten triangle
\begin{equation}
\label{equ:p}
  \xymatrix{
    \Sigma^{-1}\dual\!\Lambda \otimes_{\Lambda} P \ar[r] 
      & Y \ar[r] & C \ar[r] & {}
           }
\end{equation}
in  $\sK(\Inj\:\Lambda)$, where $P$ is a projective resolution of the
module $M$.

It is easy to see that $\Sigma^{-1}\dual\!\Lambda \otimes_{\Lambda} P$
is in $\cat$, and so \eqref{equ:p} is in particular an
Auslander-Reiten triangle in $\cat$.  As $\subcat$ is covering in
$\cat$ by Lemma \ref{lem:mod_covering}, there is a $\subcat$-cover
\[
  \xymatrix{
    A \ar[r] & \Sigma^{-1}\dual\!\Lambda \otimes_{\Lambda} P.
           }
\]
But now Theorem \ref{thm:main} says that there is an Auslander-Reiten
triangle
\[
  \xymatrix{
    A \ar[r] & B \ar[r] & C \ar[r] & {}
           }
\]
in $\subcat$, and by Remark \ref{rmk:mod} this gives an
Auslander-Reiten sequence
\[
  \xymatrix{
    0 \ar[r] & \H^0(A) \ar[r] & \H^0(B) \ar[r] & \H^0(C) \ar[r] & 0
           }
\]
in $\mod(\Lambda)$.  Since $\H^0(C) \cong M$, this gives the first
Auslander-Reiten sequence claimed in the theorem.
\end{proof}

\medskip
\noindent
{\em Acknowledgement. } I thank Henning Krause for showing me how to
use the methods from \cite{Krause3} and \cite{KrauseLe} to remove the
assumption $\gldim \Lambda < \infty$ from Section \ref{sec:example},
and Karin Erdmann and Vanessa Miemitz for their interest in this work.

\end{document}